\documentclass[12pt, titlepage]{amsart}

\usepackage{ae} % or {zefonts}
\usepackage[T1]{fontenc}
\usepackage[cp1250]{inputenc}
\usepackage{amsmath}
\usepackage{amssymb, amsfonts,amscd,verbatim}
\usepackage[dvips]{graphicx}
\usepackage{latexsym}
\usepackage{indentfirst}
\usepackage{latexsym}

\usepackage{graphicx}
\usepackage{verbatim}   % useful for program listings
\usepackage{color}      % use if color is used in text

\usepackage{amsmath}    % need for subequations

\theoremstyle{plain}
\newtheorem{Prop}{Proposition}[section]
\newtheorem{Thm}[Prop]{Theorem}
\newtheorem{Cor}[Prop]{Corollary}

\newtheorem{Lem}[Prop]{Lemma}

\theoremstyle{definition}
\newtheorem{Def}[Prop]{Definition}

\theoremstyle{remark}

\newtheorem{Problem}[Prop]{\bf Problem}
\newtheorem{Question}[Prop]{\bf Question}

\def\int{\mathop{\roman{int}}}

\def\1{^{-1}}

\def\Z{{\mathbf Z}}

\def\R{{\mathbf R}}

\def\dim{\text{dim}}

\def\dist{\text{dist}}

\def\asdim{\mathrm{asdim}}

\def\dim{\mathrm{dim}}

\def\dokaz{{\bf Proof. }}
\def\edokaz{\hfill $\blacksquare$}

\numberwithin{equation}{section}
\begin{document}
%\fontsize{20}{34pt}\selectfont
\title[
Assouad-Nagata dimension of wreath products of groups
]%
   {Assouad-Nagata dimension of wreath products of groups}

\date{ \today
}

\author{N.~Brodskiy}
\address{University of Tennessee, Knoxville, TN 37996, USA}
\email{brodskiy@math.utk.edu}

\author{J.~Dydak}
\address{University of Tennessee, Knoxville, TN 37996, USA}
\email{dydak@math.utk.edu}
\author{U.~Lang}
\address{Eidgen Technische Hochschule Zentrum, CH-8092 Z\" urich, Switzerland}
\email{lang@math.ethz.ch}

\thanks{The second-named author was partially supported
by the Center for Advanced Studies in Mathematics
at Ben Gurion University of the Negev (Beer-Sheva, Israel).}

\begin{abstract}
Consider the wreath product $H\wr G$, where $H\ne 1$ is finite and $G$ is finitely generated.
We show that the Assouad-Nagata dimension $\dim_{AN}(H\wr G)$ of $H\wr G$
depends on the growth of $G$ as follows:
\par If the growth of $G$ is not bounded by a linear function, then $\dim_{AN}(H\wr G)=\infty$,
otherwise $\dim_{AN}(H\wr G)=\dim_{AN}(G)\leq 1$.
\end{abstract}
\maketitle

%\medskip
%Printed on \today.
%\medskip
\tableofcontents

\section{Introduction}

 P.Nowak \cite{Nowak} proved that the Assouad-Nagata
dimension of some wreath products $H\wr G$ is infinite,
where $H$ is finite and $G$ is a finitely generated amenable group whose Folner function grows sufficiently
fast and satisfies some other conditions suitable for applying
Erschler's result \cite{Er}.  That result says the Folner function
$F(H\wr G)$ of $H\wr G$ is comparable to $F(H)^{F(G)}$
and the passage from it to Assouad-Nagata dimension of $H\wr G$
is fairly complicated as it includes Property A.
Thus, results of \cite{Nowak} apply only to amenable groups $G$
and do not apply neither to lamplighter groups (as the Folner function of $\Z$ is linear)
nor to wreath products with free non-Abelian groups (as those are not amenable).
\par In this paper we show that the Assouad-Nagata dimension of $H\wr G$ completely
depends on the linearity of the growth of $G$.
In particular, the lamplighter groups are not finitely presented and are of Assouad-Nagata dimension $1$
which solves positively the following question of \cite{DH}:
\begin{Question}
Is there a finitely generated group of Assouad-Nagata dimension
$1$ that is not finitely presented?
\end{Question}

\section{Assouad-Nagata dimension}

Let $X$ be a metric space and $n\ge 0$.
An {\it $n$-dimensional control function of $X$}
is a function $D_X^n\colon \R_+\to \R_+\cup\infty$ with the following property:
\par For any $r > 0$
there is a cover $\{X_0,\ldots, X_n\}$ of $X$
whose Lebesque number is at least $r$ (that means every open $r$- ball $B(x,r)$
is contained in some $X_i$) and every $r$-component of $X_i$
is of diameter at most $D^n_X(r)$. Two points $x$ and $y$ belong
to the same $r$-component of $X_i$ if there is
a sequence $x_0=x, x_1,\ldots, x_k=y$ in $X_i$ such that
$\dist(x_j,x_{j+1}) < r$ (such a sequence will be called an {\it $r$-path}).

The {\it asymptotic dimension } $\asdim(X)$ is the smallest
integer such that $X$ has an $n$-dimensional control function
whose values are finite.

The {\it Assouad-Nagata dimension} $\dim_{AN}(X)$ of a metric space $X$
is the smallest integer $n$ such that $X$ has an $n$-dimensional control function
that is a dilation ($D_X^n(r)=C\cdot r$ for some $C > 0$).

The {\it asymptotic Assouad-Nagata dimension} $\asdim_{AN}(X)$ of a metric space $X$
is the smallest integer $n$ such that $X$ has an $n$-dimensional control function
that is linear ($D_X^n(r)=C\cdot r+C$ for some $C > 0$).

In case of metrically discrete spaces $X$ (that means there is $\epsilon > 0$
such that every two distinct points have the distance at least $\epsilon$)
$\asdim_{AN}(X)=\dim_{AN}(X)$ (see \cite{BDHM}).
In particular, in case of finitely generated groups we can talk about
Assouad-Nagata dimension instead of asymptotic
Assouad-Nagata dimension.

A countable group $G$ is called {\it locally finite} if every finitely generated subgroup of $G$ is finite. A group $G$ has asymptotic dimension $0$ if and only if it is locally finite~\cite{Smith}.

Notice that $\dim_{AN}(X)=0$ if there is $C > 0$ such that
for every $r$-path the distance between its end-points is less than $C\cdot r$.
In case of groups one has the following useful criterion
of being $0$-dimensional:

\begin{Prop}\label{AsNagDimOfLocFinite}
Let $(G,d_G)$ be a group equipped with a proper left-invariant
metric $d_G$ (that means bounded sets are finite). If $G$ is locally finite,
then the following conditions are equivalent:
\begin{itemize}
\item[a.] $\dim_{AN}(G,d_G)=0$.
\item[b.] There is a constant $c > 0$ such that
for each $r > 0$ the subgroup of $G$ generated by $B(1,r)$
is contained in $B(1,c\cdot r)$.
\end{itemize}
\end{Prop}
\dokaz a)$\implies$b). Consider a constant $K > 0$
such that for each $r > 0$ all $r$-components of $G$
have diameter less than $K\cdot r$.
Notice that if $g\in G$ belongs to $r$-component of $1$
and $h\in B(1,r)$,
then $d_G(g,gh)=d_G(1,h) < r$, so $gh$ lies in the $r$-component
of $1$.
Therefore the subgroup generated by $B(1,r)$
is contained in $B(1,K\cdot r)$.
\par
b)$\implies$a). Let $G_r$ be the subgroup of $G$ generated by $B(1,r)$.
Consider two different left cosets $y \cdot G_r$ and $z \cdot G_r$ of $G_r$ in $G$.
If $d_G(yg,zh) < r$ for some $g,h\in G_r$,
then $f=h^{-1}z^{-1}yg\in B(1,r)\subset G_r$, so $y=z (hfg^{-1})$, a contradiction.
That means each $r$-component of $G$ is contained in a left coset of $G_r$
and its diameter is less than $2cr$, i.e. $\dim_{AN}(G,d_G)=0$.
\edokaz

Let us generalize $r$-paths as follows:
\par
By an {\it $r$-cube} in a metric space $X$
we mean a function $f\colon \{0,1,\ldots,k\}^n\to X$
with the property that the distance between $f(x)$ and $f(x+e_i)$
is less than $r$ for all $x \in \{0,1,\ldots,k\}^n$
such that $x+e_i\in \{0,1,\ldots,k\}^n$. Here $e_i$ belongs to the standard
basis of $\R^n$.

A sufficient condition for $\dim_{AN}(X)$ being positive
is the existence, for every $C > 0$, of an $r$-path joining
points of distance at least $C\cdot r$.
The purpose of the remainder of this section is to find a similar
sufficient condition for $\dim_{AN}(X)\ge n$.

\begin{Lem}\label{LatticeLemma}
Consider the set $X=\{0,1,\ldots,k\}^n$ equipped with the $l_1$-metric.
Suppose $X=X_1\cup\ldots\cup X_n$.
If the open $(n+1)$-ball of every point of $X$ is contained in some $X_i$,
then a $2$-component of some $X_i$
contains two points whose $i$-coordinates differ by $k$.
\end{Lem}
\dokaz Let us proceed by contradiction
and assume all $2$-components of each $X_i$
do not contain points whose $i$-coordinates differ by $k$.
Create the cover $A_i$, $1\leq i\leq n$, of the solid cube $[0,k]^n$
by adding unit cubes to $A_i$ whenever all vertices of it are contained
in $X_i$.
Given $i\in \{1,\ldots,n\}$ consider the two faces $L_i$ and $R_i$ of
$[0,k]^n$ consisting of points whose $i$-th coordinate
is $0$ and $k$ respectively.
Let $B_i$ be the complement of the $\frac{1}{4}$-neighborhood
of $A_i\cup L_i\cup R_i$. Notice that $B_i$ separates
between $L_i$ and $R_i$. Indeed,
if $L_i\cup R_i$ belongs to the same component of the $\frac{1}{4}$-neighborhood
of $A_i\cup L_i\cup R_i$, then one can find a $\frac{1}{2}$-path
in $A_i$ between points in $X_i$ whose $i$-coordinates differ by
$k$. Picking points in $X_i$ in the same unit cubes
as vertices of the path one gets a $2$-path in $X_i$
between points in $X_i$ whose $i$-th coordinates differ by $k$.

Now we get a contradiction as $\bigcap\limits_{i=1}^n B_i=\emptyset$
in violation of the well-known result in dimension theory about separation
(see Theorem 1.8.1 in \cite{Engel}).
\edokaz

\begin{Cor}\label{CubeCorollary}
Suppose $X$ is a metric space
with an $(n-1)$-dimensional control function $D_X^{n-1}\colon \R_+\to \R_+\cup\infty$.
For any $r$-cube $$f\colon \{0,1,\ldots,k\}^n\to X$$
there exist two points $a$ and $b$ in $\{0,1,\ldots,k\}^n$ whose $i$-th coordinates differ
by $k$ for some $i$ and $\dist(f(a),f(b))\leq D_X^{n-1}(n\cdot r)$.
\end{Cor}
\dokaz
Consider a cover $X=X_1\cup\ldots\cup X_n$ of $X$ of Lebesque number at least $n\cdot r$
such that $n\cdot r$-components of each $X_i$ are of diameter at most $D_X^{n-1}(n\cdot r)$. The cover $\{0,1,\ldots,k\}^n=f^{-1}(X_1)\cup \ldots\cup f^{-1}(X_n)$
has the property that the open $(n+1)$-ball of every point is contained in some $f^{-1}(X_i)$,
so by \ref{LatticeLemma} a $2$-component (in the $l_1$-metric) of some $f^{-1}(X_i)$
contains two points $a$ and $b$ whose $i$-coordinates differ by $k$.
Therefore $f(a)$ and $f(b)$ belong to the same $r$-component
of $X_i$ and $\dist(f(a),f(b)) \leq  D_X^{n-1}(n\cdot r)$.
\edokaz

We need an upper bound on the size of $r$-cubes $f$
in terms of dimension control functions and the Lipschitz constant
of $f^{-1}$. One should view the next result as a discrete analog of the fact
that one cannot embed $I^n$ into an $(n-1)$-dimensional topological space.
\begin{Cor}\label{CubeLipCorollary}
Suppose $X$ is a metric space
with an $(n-1)$-dimensional control function $D_X^{n-1}\colon \R_+\to \R_+\cup\infty$.
If $f\colon \{0,1,\ldots,k\}^n\to X$ is an $r$-cube, then
$k\leq D_X^{n-1}(n\cdot r)\cdot Lip(f^{-1})$.
\end{Cor}
\dokaz
By \ref{CubeCorollary} there is an index $i\leq n$ and points $a$ and $b$
whose $i$-coordinates differ by $k$ such that $\dist(f(a),f(b)) \leq  D_X^{n-1}(n\cdot r)$.
Since $k\leq \dist(a,b)\leq Lip(f^{-1})\cdot \dist (f(a),f(b))\leq D_X^{n-1}(n\cdot r)\cdot Lip(f^{-1})$,
we are done.
\edokaz

\section{Wreath products}

Let  $ A$ and $B$ be groups. Define the action  of $ B$
on the direct product $A^B$ (functions have finite support) by
 \begin{displaymath}
b f(\gamma) := f(b^{-1}\gamma),
\end{displaymath}
for any
 $f\in A^B$
 and  $\gamma\in B$. The {\it wreath product} of $ A$ and $ B$, denoted  $A\wr B$, is the  semidirect product $A^B\rtimes B$ of groups $A^B$ and $B$.
 That means it consists of ordered pairs $(f,b)\in A^B\times B$
 and $(f_1,b_1)\cdot (f_2,b_2)=(f_1(b_1 f_2),b_1b_2)$.

 \par We will identify $(1,b)$ with $b\in B$ and $(f_a,1)$ with $a\in A$,
 where $f_a$ is the function sending $1\in B$ to $a$ and $B\setminus \{1\}$ to $1$.
 This way both $A$ and $B$ are subgroups of $A\wr B$
 and it is generated by $B$ and elements of the form $b\cdot a\cdot b^{-1}$.
 That way the union of generating sets of $A$ and $B$ generates $A\wr B$.

The {\it lamplighter group} $L_n$ is the wreath product $\Z/n\wr \Z$ of $\Z/n$ and $\Z$.

If $g\in G$ and $a\in H\setminus \{1\}$, then $g\cdot a\cdot g^{-1}\in K$ will be called
the {\it $a$-bulb indexed by $g$}
or the {\it $(g,a)$-bulb}. A {\it bulb} is a $(g,a)$-bulb for some $a\in H$ and some $g\in G$.

Consider the wreath product $H\wr G$, where $H$ is finite and $G$ is finitely generated.
Let $K$ be the kernel of $H\wr G\to G$. $K$
is a locally finite group (the direct product of $|G|$ copies of $H$).
In case $H$ is finite we choose as a set of generators of $H\wr G$ the union of $H\setminus \{1\}$
and a set of generators of $G$.

\begin{Lem}\label{LowerBoundLem}
Suppose $n > 1$.
Any product of bulbs indexed by mutually different
elements $g_i\in G$, $i\in\{1,\ldots, n\}$, has length at least $n$.
\end{Lem}
\dokaz
Consider $x=(g_1a_1g_1^{-1})\cdot\ldots\cdot (g_na_ng_n^{-1})\in K$. If its length is smaller than $n$, then $x=x_1\cdot b_1\cdot x_2\cdot b_2\cdot\ldots\cdot x_k\cdot b_k\cdot x_{k+1}$, where $k < n$
and $b_i\in H$, $x_i\in G$ for all $i$.
We can rewrite $x$ as $(y_1\cdot b_1\cdot y_1^{-1})\cdot (y_2\cdot b_2\cdot y_2^{-1})
\cdot\ldots \cdot (y_k \cdot b_k\cdot y^{-1}_k)\cdot y$,
where $y_1=x_1$. Since $x\in K$, $y=1$.
Now we arrive at a contradiction by looking at projections of $K$ onto its summands.
\edokaz

\begin{Lem}\label{RepresentationLem}
Suppose $r > 1$.
Any element of $K$ of length less than $r$
is a product of bulbs indexed by elements of $G$ of length less than $r$.
\end{Lem}
\dokaz Any element of $K$ of length less than $r$
has $x_1a_1x_2a_2\ldots x_ka_kz$ as a minimal representation, so
it can be rewritten as
$$(x_1a_1x_1^{-1})(x_1x_2a_2x_2^{-1}x_1^{-1})\ldots
(x_1\ldots x_ka_kx_k^{-1}\ldots x_1^{-1}).$$
 Therefore the bulbs involved
are indexed by elements of $G$ of length less than $r$.
\edokaz

In case of the lamplighter group $L_2$ there is a precise calculation
of length of its elements in \cite{DeadEnds}.
We need a generalization of those calculations.
\begin{Lem}\label{UpperBoundLem}
Let $H$ be finite.
Suppose the subgroup $\Z$ generated by $t\in G$
is of finite index $n$ and there are generators $\{t,g_1,\ldots,g_n\}$ of $G$
such that every element $g$ of $G$ can be expressed as $g_i\cdot t^{e(g)}$ for some $i$.
\begin{enumerate}
\item Every element of $K$ can be expressed as a product
of $(h_i,a_i)$-bulbs, $i=1,\ldots,k$, such that $h_i\ne h_j$ for $i\ne j$.
\item The length of such product is at most $n(k+2+4\max\{|e(h_i)|\})$.
\end{enumerate}

\end{Lem}
\dokaz
Observe the product of the $(g,a)$-bulb and the $(g,b)$-bulb
is the $(g,a\cdot b)$-bulb, so every product of bulbs
can be represented as a product
of $(h_i,a_i)$-bulbs, $i=1,\ldots,k$, such that $h_i\ne h_j$ for $i\ne j$.
We will divide those bulbs in groups determined by
$h_i\cdot t^{-e(h_i)}$.
Since there are at most $n$ groups, it suffices to show that if $h_i\cdot t^{-e(h_i)}=g$ for all $i$,
then the length of the product $x$ of $(h_i,a_i)$-bulbs is at most $k+2+4\max\{|e(h_i)|\}$.
We may order $h_i$ so that the function $i\to e(h_i)$ is strictly increasing.
Now, $$g^{-1}\cdot x\cdot g=\prod\limits_{i=1}^k t^{e(h_i)}\cdot a_i\cdot t^{-e(h_i)}
=t^{e(h_1)}\cdot a_1\cdot t^{-e(h_1)+e(h_2)}\cdot a_2\cdot\ldots\cdot a_k\cdot t^{-e(h_k)}$$
and its length is at most $k+|e(h_1)|+e(h_k)-e(h_1)+|e(h_k)|\leq k+4\max\{|e(h_i)|\}$.
Therefore the length of $x$ is at most $k+2+4\max\{|e(h_i)|\}$.
\edokaz

\section{Dimension control functions of wreath products}

Recall that the {\it growth} $\gamma$ of $G$ is the function counting
the number of points in the open ball $B(1,r)$ of $G$
for all $r > 0$. Notice that $\gamma$ being bounded by a linear function
is independent on the choice of generators of $G$.

The next result relates the growth
function of $G$ to dimension control functions of the kernel of the projection $H\wr G\to G$.

\begin{Thm}\label{ANDimOfKernelB}
Suppose $G$ and $H$ are finitely generated and $K$
is the kernel of the projection $H\wr G\to G$ equipped with the metric induced from $H\wr G$.
If $\gamma$ is the growth function of $G$
and  $D_K^{n-1}$ is an $(n-1)$-dimensional control
function of $K$,
then the integer part of
$\frac{\gamma(r)}{n}$ is at most $D_K^{n-1}(3nr)$.
\end{Thm}
\dokaz
Given $k \ge 1$ we will construct
a $3r$-cube $f\colon \{0,k\}^n\to K$ similarly to the way paths
in the Cayley graph of $K$ are constructed.
There, it suffices to label the beginning vertex and all the edges
and that induces labeling of all the vertices.
In case of our $3r$-cube we label the origin by $1\in K$
and each edge from $x$ to $x+e_i$, $e_i$ being an element of the standard
basis of $\R^n$, will be labeled by $x(j,i)$, where $j$ is the $i$-th coordinate of $x$.
It remains to choose $x(j,i)$, $1\leq i\leq n$ and $0\leq j\leq k-1$.
Given $r > 0$ consider mutually different elements $g(j,i)$, $1\leq i\leq n$ and $0\leq j\leq k-1$ of $G$ whose length is smaller than $r$, where $k$ is the integer part of
$\frac{\gamma(r)}{n}$.
Pick $u\in H\setminus \{1\}$ and put $x(j,i)=g(j,i)\cdot u\cdot g(j,i)^{-1}$.
By \ref{LowerBoundLem} one has $Lip(f^{-1})\leq 1$,
so $k \leq D_K^{n-1}(3nr)$ by \ref{CubeLipCorollary}.
\edokaz

If $H$ is finite, then the kernel $K$ of the projection $H\wr G\to G$
is locally finite and it has a $0$-dimensional control
function $D^0_K$ attaining finite values ($K$ is equipped
with the metric induced from $H\wr G$).
Let us relate $D^0_K$ to the growth of $G$.

\begin{Thm}\label{0DimControlOfK}
Suppose $G$ is finitely generated and $H\ne\{1\}$ is finite. Let $K$
be the kernel of the projection $H\wr G\to G$ equipped with the metric induced from $H\wr G$.
If $\gamma$ is the growth function of $G$,
then $D_K^{0}(r):=(2r+1)\gamma(r)$ is a $0$-dimensional control
function of $K$.
\end{Thm}
\dokaz
It suffices to show that $r$-component of $1$ in $K$ is of diameter
at most $(2r+1)\gamma(r)$ as any $r$-component of $K$ is a shift
of the $r$-component containing $1$.
By \ref{RepresentationLem} any element of $B(1,r)$ in $K$ is a product
of bulbs indexed by elements of $G$ of length less than $r$.
Therefore any product of elements in $B(1,r)$
is a product of bulbs indexed by elements of $G$ of length less than $r$
and such product can be reduced to a product
of at most $\gamma(r)$ such bulbs.
Each of them is of length at most $2r+1$, so the length of the product
is at most $(2r+1)\cdot\gamma(r)$.
\edokaz

\begin{Thm}[cf. {\cite[Proposition 4.2]{Dra}}]\label{nDimControlOfKernel}
Suppose $G$ is finitely generated and $\pi\colon G\to I$ is a
retraction onto its subgroup $I$ with kernel $K$. $K$ is equipped
with the metric induced from a word metric on $G$ so that
generators of $I$ are included in the set of generators of $G$. If
$D^n_I$ is an $n$-dimensional control function of $I$ and $D^0_K$
is a $0$-dimensional control function of $K$, then
$$D^n_I(r)+D^0_K(r+2D^n_I(r))$$
 is an $n$-dimensional control
function of $G$.
\end{Thm}
\dokaz
Given $r > 0$ express $I$ as $I_0\cup\ldots\cup I_n$
so that $r$-components of $I_i$ have diameter at most $D^n_I(r)$.
Consider $G_i=\pi^{-1}(I_i)$.
If $g_1\cdot 1,\ldots,g_1\cdot x_m$ is an $r$-path in $G_i$,
then $h_1=\pi(g_1)\cdot 1,\ldots,h_m=\pi(g_1)\cdot y_m$ form an $r$-path in $I_i$
(here $y_j=\pi(x_j)$),
so $l(y_j)\leq D^n_I(r)$ for all $j$.
Consider $z_j=x_j\cdot y_j^{-1}\in K$.
Notice $\dist(z_j,z_{j+1})< r+2D^n_I(r)$.
Therefore, $\dist(1,z_m)\leq D^0_K(r+2D^n_I(r))$
resulting in $l(x_m) \leq D^0_K(r+2D^n_I(r))+D^n_I(r)$
and $\dist(g_1,g_1\cdot x_m) \leq D^n_I(r)+D^0_K(r+2D^n_I(r))$
which completes the proof.
\edokaz

\begin{Def}[cf. {\cite[Section VI.B]{Harpe}}]
Let $f$ and $g$ be functions from $\R_+$ to $\R_+$. We say that
$f$ {\it weakly dominates} $g$ if there exist constants
$\lambda\ge 1$ and $C\ge 0$ such that $g(t)\le \lambda f(\lambda
t+C)+C$ for all $t\in \R_+$.

Two functions are {\it weakly equivalent} if each weakly dominates
the other.
\end{Def}

\begin{Thm}
Suppose $G$ is finitely generated infinite group and $H\ne\{1\}$
is finite. Let $\gamma$ be the growth function of G and $D^n_G$ be
an $n$-dimensional control function of $G$. Then for any $k\ge n$
there is a $k$-dimensional control function of $H\wr G$ which is
weakly dominated by $(D^n_G(t)+t)\cdot\gamma(D^n_G(t)+t)$. Also,
for any $k\ge n$ every $k$-dimensional control function of $H\wr
G$ weakly dominates the function $\gamma$.
\end{Thm}

\dokaz Notice that $\gamma$ dominates a linear function and
combine~\ref{0DimControlOfK} and~\ref{nDimControlOfKernel}. To get
the estimate from below, notice that a $k$-dimensional control
function of $H\wr G$ works as a $k$-dimensional control function
of the kernel $K$, and apply~\ref{ANDimOfKernelB}.
\edokaz

Our next result gives a better solution to Question 2
in~\cite{Nowak}.

\begin{Cor}
Suppose $G$ is a finitely generated group of exponential growth
and $H\ne\{1\}$ is finite. If $\dim_{AN}(G)\le n$ then for any
$k\ge n$ the $k$-dimensional control function of $H\wr G$ is
weakly equivalent to the function $2^t$ (i.e. there is a
$k$-dimensional control function of $H\wr G$ weakly dominated by
$2^t$ and every such control function weakly dominates $2^t$).
\end{Cor}

\begin{Cor}\label{nDimControlOfWreathWithFree}
Let $F_2$ be the free non-Abelian group of two generators. For
every $n\ge 1$ the $n$-dimensional control function of $\Z/2\wr
F_2$ is weakly equivalent to the function $2^t$ (i.e. there is an
$n$-dimensional control function of $\Z/2\wr F_2$ weakly dominated
by $2^t$ and every such control function weakly dominates $2^t$).
\end{Cor}

\dokaz Notice that the function $f(t)=2^t$ is weakly equivalent to
the growth function of $F_2$ and $\dim_{AN}(F_2)=1$. \edokaz

\section{Assouad-Nagata dimension of wreath products}

Suppose $G$ is finitely generated and $H\ne 1$ is finite.
If $\dim_{AN}(G)=0$, then $G$ is finite and so is $H\wr G$.
In such case $\dim_{AN}(H\wr G)=0=\dim_{AN}(G)$.
Therefore it remains to consider the case of infinite groups $G$.

\begin{Thm}\label{LinearGrowthCase}
Suppose $G$ is an infinite finitely generated group and $H$ is a finite group.
If the growth of $G$ is bounded by a linear function, then $\dim_{AN}(K)=0$
and $\dim_{AN}(H\wr G)=\dim_{AN}(G)=1$.
\end{Thm}
\dokaz
Notice that \ref{0DimControlOfK} does provide
a $0$-dimensional control function for $K$. However, it may not be bounded by a linear
function,
so we have to do more precise calculations.
\par
$G$ is a virtually nilpotent group by Gromov's Theorem (see \cite{GroNilp} or Theorem 97 in \cite{Kap}).
Let $F$ be a nilpotent subgroup of $G$ of finite index.
Pick elements $a_i$, $i=1,\ldots,k$, of $G$ such that
$G=\bigcup\limits_{i=1}^k a_i\cdot F$ and pick a natural $n$ satisfying $|a_i|\leq n$ for all $i\leq k$.
Every two elements of $F$ can be connected in $G$ by a $2$-path.
From each point of the path (other than initial and terminal points)
one can move to $F$ by distance at most $n$ (by representing that
point as $a_i\cdot x$ for some $x\in F$). Therefore we can create
a $(2n+2)$-path in $F$ joining the original points.
That means $F$ is generated by its elements of length at most $2n+1$,
\par
Let $\{F_i\}$ be the lower central series of $F$
and let $d_i$ be the rank of $F_i/F_{i+1}$, $i\ge 0$.
Since the growth of $F$
is also linear, Bass' Theorem (see \cite{Bass} or Theorem 103 in \cite{Kap})
stating that the growth of $F$ is polynomial of degree
$d=\sum\limits_{i=0}^\infty (i+1)\cdot d_i$
implies that $d_0=1$ and all the other ranks $d_i$ are $0$.
Hence the abelianization of $F$ is of the form
$\Z\times A$, $A$ being a finite group, and the commutator
group of $F$ is finite.
Therefore $F$ is virtually $\Z$ and that means $G$ is virtually $\Z$ as well.
\par
Let $n$ be the index of $\Z$ in $G$
and pick elements $g_1,\ldots, g_n$ of $G$
such that any element of $G$ can be expressed
as $g_i\cdot t^k$ for some $i\leq n$ and some $k$, where $t$ is the generator of $\Z\subset G$.
Without loss of generality we may assume that
the set of generators of $G$ chosen to compute the word length $l(w)$
of elements $w\in H\wr G$ is $t,g_1,\ldots, g_n$.
For $H$ we choose all of $H\setminus \{1\}$ as the set of generators.

We need existence of $C > 0$ such that
 $\frac{|k|}{C}\leq l(t^k)\leq |k|$ for all $k$.
 It suffices to consider $k > 0$.
 Since the number of points in $B(1_G,4)$ is finite,
 there is $C > 0$ such that $t^u\in B(1_G,4)$ implies $|u| \leq C$.
Now, if $l(t^k)=m$ and $t^k=x_1\cdot\ldots \cdot x_m$,
where $l(x_i)=1$, then there are $u(i)$ such that
$\dist(x_1\cdot\ldots\cdot x_i,t^{u(i)})\leq 1$ for all $i\leq k$
(we choose $u(m)=k$ obviously).
Therefore $\dist(t^{u(i)},t^{u(i+1)})\leq 3$
and $u(i+1)-u(i)\leq C$.
Now $k=u(m)=(u(m)-u(m-1))+\ldots+(u(2)-u(1))+u(1)\leq C\cdot m$
implying $l(t^k)=m\ge \frac{k}{C}$.

\par
By \ref{RepresentationLem} any element of $K$ of length less than $r$
is a product of bulbs
indexed by elements of $G$ of length less than $r > 1$.
If $l(g_i\cdot t^k) < r$, then $l(t^k)< r+1< 2r$
and $|k|\leq C\cdot l(t^k)\leq 2Cr$.
Therefore there are at most $n\cdot 4Cr$ such words
and
any product of such bulbs is of length at most
$n(4Crn+2+2Cr)\leq r(4Cn^2+2n+2Cn)$ by \ref{UpperBoundLem}.

Therefore the group
generated by $B(1,r)$ in $K$ is contained in $B(1,Lr)$,
where $L=4Cn^2+2n+2Cn$,
and $\dim_{AN}(K)=0$ by \ref{AsNagDimOfLocFinite}.
Using the Hurewicz Theorem for Assouad-Nagata dimension of \cite{BDLM}
we get $\dim_{AN}(H\wr G)\leq \dim_{AN}(G)=1$ (one can also use
\ref{nDimControlOfKernel}).
Since $H\wr G$ is infinite, its Assouad-Nagata dimension is positive
and $\dim_{AN}(H\wr G)=\dim_{AN}(G)=1$.
\edokaz

\begin{Cor}\label{NonLinearGrowth}
If the growth of $G$ is not bounded by a linear function and $H\ne 1$, then $\dim_{AN}(H\wr G)=\infty$.
\end{Cor}
\dokaz
Let $\gamma$ be the growth of $G$ in some set of generators.
Suppose $\dim_{AN}(K) <n < \infty$,
so it has an $(n-1)$-dimensional function
of the form $D_K^{n-1}(r)=C\cdot r$ for some $C > 0$.
By \ref{ANDimOfKernelB} one has
$\gamma(r)/n\leq C\cdot 3nr +1$.
Thus $\gamma(r)\leq n\cdot (3nCr+1)$ and the growth of $G$ is bounded by a linear
function,
a contradiction.
\edokaz

\begin{Problem}
Suppose $G$ is a locally finite group equipped with a proper
left-invariant metric $d_G$. If $\dim_{AN}(G,d_G) > 0$,
is $\dim_{AN}(G,d_G)$ infinite?
\end{Problem}

\end{document}